\documentclass[12pt]{article}
\usepackage{amsfonts,amssymb,amscd,amsthm}

\newcommand{\sect}[1]{\setcounter{equation}{0}\section{#1}}

\def\be{\begin{equation}}
\def\ee{\end{equation}}
\def\bea{\begin{eqnarray}}
\def\eea{\end{eqnarray}}

\def\C{{\mathbb C}}

\def\hbar{{\xi}}
\def\z{{z}}

\def\conm#1#2{\left [ {#1},{#2} \right ]}

\def\anticonm#1#2{\left\{ {#1},{#2} \right\} }

\begin{document}


\thispagestyle{empty}

\hfill \today 

\vspace{2.5cm}

\begin{center}
\sc{\LARGE From Quantum Universal Enveloping Algebras to Quantum Algebras}
\end{center}

\bigskip\bigskip

\begin{center}
E. Celeghini$^1$,  A. Ballesteros$^2$ and M.A. del Olmo$^3$
\end{center}

\begin{center}
$^1${\sl Departimento di Fisica, Universit\`a  di Firenze and
INFN--Sezione di
Firenze \\
I50019 Sesto Fiorentino,  Firenze, Italy}\\
\medskip

$^2${\sl Departamento de F\'{\i}sica, Universidad de Burgos, \\
E-09006, Burgos, Spain.}\\
\medskip

$^3${\sl Departamento de F\'{\i}sica Te\'orica, Universidad de
Valladolid, \\
E-47005, Valladolid, Spain.}\\
\medskip

{e-mail: celeghini@fi.infn.it, angelb@ubu.es, olmo@fta.uva.es}

\end{center}

\bigskip

\bigskip

\begin{abstract}

The ``local'' structure of a quantum group $G_q$ is currently
considered to be an infinite-dimensional object: the corresponding
quantum universal enveloping algebra $U_q(g)$, which is a Hopf
algebra deformation of the universal enveloping algebra  of a
$n$-dimensional Lie algebra $g=Lie(G)$. However, we show how, by
starting from the generators of the underlying Lie bialgebra
$(g,\delta)$, the analyticity in the deformation parameter(s)
allows us to determine in a unique way a set of $n$ ``almost
primitive'' basic objects in $U_q(g)$, that could be properly
called the ``quantum algebra generators''. So, the analytical
prolongation $(g_{q},\Delta)$ of the Lie bialgebra $(g,\delta)$ is
proposed as the appropriate local structure of $G_q$.
Besides, as in this way  $(g,\delta)$ and $U_q(g)$ are shown to be
in one-to-one correspondence,
  the classification of quantum groups is reduced to the
classification of Lie bialgebras. The $su_q(2)$ and $su_q(3)$ cases are explicitly
elaborated.

\end{abstract}

\vskip 1cm


\vskip 0.4cm


\vfill\eject


\sect{Introduction}

Quantum groups are a non-commutative generalization of Lie groups
endowed with a Hopf algebra structure~\cite{CP,fuchs}. Some attempts in
order to get structural properties of these objects  have been
previously considered (see \cite{LM,mudrov} for a prescription to
get the quantum coproduct --but not the deformed commutation
rules-- for a wide class of examples). Moreover, to our knowledge,
a general investigation concerning the uniqueness of this
quantization process has not been given yet and only restrictive
results for certain deformations of simple Lie algebras have been
obtained (see \cite{schnider}, Chapter 11). As a consequence of
the above mentioned facts, a complete classification of quantum
groups in the spirit of Cartan cannot be found in the literature
(see \cite{ballesteros04,ballesteros05} and references therein).

We would like to stress that both in Lie group theory as well as
in their physical applications, the infinitesimal counterpart of a
Lie group transformation --i.e. its Lie algebra $g$-- plays a
fundamental role, since it is interpreted as the local (around the
identity) symmetry. Correspondingly, the local counterpart of
quantum groups was soon algebraically identified through Hopf
algebra duality, presents a wealth of interesting mathematical
properties, and has been also applied in different physical
contexts.

However, such infinitesimal counterpart of a quantum group $G_q$
is not a deformed Lie algebra, but a quantum universal enveloping
algebra $U_q(g)$: a Hopf algebra deformation of the universal
enveloping algebra $U(g)$ of $g$, i.e., a deformation of the
infinite-dimensional object that has as a basis the set of all
ordered monomials of powers of the generators of $g$. This means that
quantum deformations lead us locally  (therefore, geometrically)
to structures which are quite different from Lie algebras. In
particular, when we consider a deformation of an $n$ dimensional
Lie algebra, only the infinite dimensional algebra $U_q(g)$ makes sense in
despite that its Poincar\'e-Birkhoff-Witt (PBW) basis is
constructed in terms of a basic set of elements of the same
dimension of the Lie algebra. Indeed, contrarily to the non
deformed case, where inside all sets of basic elements the vector
space of generators is univocally defined, in $U_q(g)$ there is an
infinitude of basic sets coexisting on the same footing.

 This problem of the basis underlies many difficulties encountered
when a precise physical/geometrical meaning has to be assigned to
the $U_q(g)$, as in the context of quantum deformations of
space-time symmetries. In that case, it is well known that the
models so obtained depend on the choice of different bases
(for instance, the bicrossproduct one \cite{majid}), and that different possibilities
are related through non-linear transformations. The aim of this
paper is to solve this problem providing a universal and computational prescription
for the characterization and construction of the $n$-dimensional
quantum analogue of a Lie Algebra.

To do this, we have to analyse the role and properties of the Lie
algebra generators  within $U(g)$. Among the infinite possible PBW
bases, all of them related by non-linear invertible transformations, the
generators determine the only one closed under linear commutation
rules and whose tensor product representations are constructed
additively. The latter property can be stated in Hopf-algebraic
terms as the Friedrichs theorem~\cite{jacobson}: the only
primitive elements $\{X_j\}$ in $U(g)$ ({\it i.e.} the elements such
that $\Delta(X_j)=\Delta_{(0)}(X_j):= 1\otimes X_j + X_j\otimes 1$)
are just the generators of $g$ as a Lie algebra. In this way, the
generators of $g$ become distinguished elements of $U(g)$.

The  additivity of generators in   representation theory  is the
reason why in physical applications we are used to disregard
$U(g)$ as a mathematical curiosity and to focus on the quite more
manageable Lie algebra. However, we realize immediately that the
situation changes drastically in $U_q(g)$, where the law for the
construction of tensor product representations (coproduct)
includes non-linear functions and no primitive bases exist. In
this paper we show that, among the infinitely many possible bases,
there is only one (that we will call ``almost primitive" basis)  where the
coproducts are ``as primitive as possible'', since all inessential
terms have been removed. Indeed, the only changes from the
primitive coproducts are those imposed by the bialgebra
cocommutator $\delta$ to be consistent with the Hopf algebra
postulates. This almost primitive basis is proposed as the true
deformation of the Lie algebra and, thus, called {\it quantum
algebra}.

It is well-known that Lie
group theory is based on analyticity with respect to group
parameters. In the same way, analyticity in the deformation
parameter(s) will give us the keystone for the identification of
the proper quantum algebra, that will be defined as the $n$-{\it
dimensional vector space $(g_q,\Delta) \subset U_q(g)$ obtained as
analytical prolongation of the Lie bialgebra $(g,\delta)$} (note
that  analyticity in the deformation has played yet a useful role
in quantum algebras, for instance in their contractions
\cite{LED}).

In this analytical prolongation, the Lie bialgebra cocommutator
map \,$\delta$ describes the first order deformation and can be
considered as the derivative at the origin of the quantum
coproduct. This $\delta$, together  with the zero-order
deformation (the Lie-Hopf algebra) and the coassociativity of the
coproduct will allow us to construct order-by-order the deformed
coproduct. The commutators (by inspection, not $q$-commutators)
are then obtained imposing the homomorphism for the coproduct.

Summarizing, in this paper we are attempting to describe the
``commutative'' (in a broad sense) diagram
\[ \begin{CD}
(g,\Delta_{(0)}) @<{\rm \hskip0.42cm Friedrichs\; theorem\hskip0.42cm}<< U(g)  \\
@VV{q}V @ VV{q}V  \\ (g_q,\Delta) @<{\rm Generalized\;
Friedrichs\; th.}<< U_q(g)
\end{CD}\]
\noindent where a new object, the quantum algebra $(g_q,\Delta)$,
is introduced, and its connections with its neighbours in the
diagram fully discussed. The vertical lines of the diagram
represent the quantization procedure, and the horizontal ones are
related to the definition of the basic set of the universal
enveloping algebras and their quantum analogues. Remember that,
for a given Lie-Hopf algebra $(g,\Delta_{(0)})$, several Lie
bialgebras $(g,\delta)$ exist and each of them determines one
different quantization and, as a consequence, a different diagram.

The paper is organized as follows.
In section~\ref{analyticalquantization} we describe the
analytical approach to the problem and, in particular, the
relation between the Lie algebra $(g,\Delta_{(0)})$ and its
analytical prolongation in the direction of $\delta$, the quantum
algebra $(g_q,\Delta)$. Moreover, in order to make the approach
more clear, we present in Section~\ref{quantizationofsu(2)} the construction of the
standard deformation of $su(2)$. Section~\ref{quantizationofu(3)} is a true application
as it exhibits the standard quantization of $su(3)$ with all
generators, commutation relations and related coproducts.
Section~\ref{Friedrichstheoremrevisited} is devoted to revisit the first horizontal line of the diagram,
i.e. the one-to-one connection between $U(g)$ and
$(g,\Delta_{(0)})$ (Friedrichs theorem) in such a way that it can
be generalized to connect --always in a one-to-one way-- $U_q(g)$
and $(g_q,\Delta)$, a subject that is discussed in
Section~\ref{ExtensionFriedrichstheorem}.
Finally,  some  conclusions close the paper.

\sect{Analytical quantization: $(g,\Delta_{(0)}) \rightarrow
(g_q,\Delta)$}\label{analyticalquantization}

As it is well known, the quantum universal enveloping algebra
$U_q(g)$ is a Hopf algebra that depends on one deformation
parameter $z=\log q$\, (the generalization to multiparametric
deformations is straightforward) and such that, in the limit $z\to
0$ (or $q\to 1$), $U_q(g)$ becomes  $U(g)$ and all possible sets
of basic elements reduce to a basis in $g$. Also, $U_q(g)$ is the quantization of a given Lie bialgebra $(g,\delta)$ where $g$
is a Lie algebra of dimension $n$, and $\delta: g\rightarrow
g\otimes g$ is a compatible skew-symmetric map~\cite{CP}. In
particular, $U_q(g)$ is a Hopf algebra such that
\[
\delta =\lim_{z\to 0}\frac {\Delta -\sigma\circ \Delta}{2 z} ,
\]
$\sigma$ being the flip operator (i.e., $\sigma (A\otimes
B)=B\otimes A$). So, $\delta$ can be interpreted as the derivative
at the origin of the quantization and $U_q(g)$ is sometimes called
a ``quantization of $U(g)$ in the direction of $\delta$''.
Such quantization is usually constructed starting
from any PBW basis in $U(g)$. Thus, a univocal correspondence between $U_q(g)$ and
$(g,\delta)$ is found, while no general results concerning the uniqueness of the quantization process are known.

Here we present a different three-step quantization procedure:
\begin{enumerate}
\item
By using analyticity and coassociativity we find order-by-order
the changes induced in $\Delta_{(0)}$ by $\delta$ and we determine
in this way the full quantum coproduct $\Delta$.
\item
 By using analyticity and the homomorphism property of $\Delta$, we obtain the
 commutation rules for $g_q$ starting from the known ones of $g$. Thus, the
 $n$-dimensional $(g_q,\Delta)$ is constructed.
\item
A PBW basis in $U_q(g)$ is built from $g_q$.
\end{enumerate}

We thus construct an unique correspondence $(g,\delta) \rightarrow
(g_q,\Delta) \rightarrow U_q(g)$. Since $(g,\delta)$ is the limit
of $(g_q,\Delta)$ and, as shown in
Sect.~\ref{ExtensionFriedrichstheorem}, $g_q$ is the only one
almost primitive basis in $U_q(g)$ (exactly like $g$ is the only
primitive basis in $U(g)$), a one-to-one correspondence is found
between $(g,\delta)$ and $U_q(g)$. So equivalences in
$U_q(g)$ implies equivalences in the Lie bialgebras, and the
classification of $U_q(g)$ is carried to the quite simpler
classification of Lie bialgebras.

The two main assumptions of the  analytical quantization
procedure are:
\begin{enumerate}
\item
The commutation relations of any basic set $\{Y_j\}$
$(j=1,2,\dots, n)$ of $U_q(g)$ (as well as of $U(g)$) are
analytical functions of the $Y_j$.
\item
The quantum coproduct $\Delta$ of the $Y_j$ can be written as a formal
series
\be\label{coproductseries}
\Delta(Y_i)=
\sum_{k=0}^{\infty}\Delta_{(k)}(Y_i)= \Delta_{(0)}(Y_i) +
\Delta_{(1)}(Y_i)  + \dots
\ee
with \,$\Delta_{(k)}(Y_i)$\,
a homogeneous polynomial of degree $k+1$ in \,$1\otimes Y_j$\, and
\,$Y_j \otimes 1$.
\end{enumerate}

Since we are deal with a Hopf algebra, $\Delta$ has to verify the
coassociativity condition \be \label{coass} (\Delta\otimes
1-1\otimes\Delta )\circ\Delta (Y_i)=0 \ee as well as the
homomorphism property \be\label{homomorphismproperty1}
\Delta(\conm{Y_i}{Y_j})= \conm{\Delta(Y_i)}{\Delta(Y_j)} . \ee
Taking into account eq. (\ref{coproductseries}), eqs. (\ref{coass})
and (\ref{homomorphismproperty1}) can be rewritten as
\be\label{recursionk} \sum_{j=0}^{k}\left( \Delta_{(j)}\otimes 1 -
1 \otimes \Delta_{(j)} \right) \circ \Delta_{(k-j)} (Y_i)=0
,\qquad\qquad \forall k , \ee
 \be\label{homomorphismproperty}
\Delta_{(k)}(\conm{Y_i}{Y_j})= \sum_{l=0}^k
\conm{\Delta_{(l)}(Y_i)}{\Delta_{(k-l)}(Y_j)} ,\qquad\qquad
\forall k \,. \ee Note that commutation rules and coproducts in
$U_q(g)$ are fully defined from commutation rules and coproducts
in $\{Y_j\}$.

In order to deform $(g,\Delta_{(0)}) \rightarrow (g_q,\Delta)$ in
the direction of $\delta$ we have to introduce the modifications
to $\Delta_{(0)}$ imposed by the Lie bialgebra $(g,\delta)$ to be
consistent with the coassociativity. Thus we define
$$
\Delta_{(1)}(X_i) := z\; \delta (X_i),
$$
by putting to zero arbitrary
cocommutative contributions to $\Delta_{(1)}(X_i)$  because they are unrelated to $\delta$
 (see
Sect.~\ref{Friedrichstheoremrevisited} and
\ref{ExtensionFriedrichstheorem}).
Then we can write
\be\label{coproductseriesdelta}
 \Delta(X_i)=
\Delta_{(0)}(X_i)+ z\; \delta (X_i) + {\cal O}_{(2)}(X_i), \ee
where ${\cal O}_{(m)}(X_i)$ is a series of degree greater than $m$ in $
X_j\otimes 1$ and $1\otimes X_j$. Because of
(\ref{coproductseries}), ${\cal O}_{(2)}(X_i)$ can be written
\be\label{O2} {\cal O}_{(2)}(X_i)= \Delta_{(2)}(X_i)+ {\cal
O}_{(3)}(X_i) . \ee By consistency  with the coassociativity
condition  (\ref{recursionk})
for $k=2$, $\Delta_{(2)}(X_i)$ contributions must satisfy a set of
well precise conditions.  The contribution in $z^2$ is determined
by $\delta$, while the other possible contributions are consistent
with zero because they are proportional to arbitrary parameters which are independent from $\delta$, as
described in Sect.~\ref{ExtensionFriedrichstheorem}. As the
analytical procedure requires to include only the changes imposed
by $\delta$, all these last contributions are put to zero. Thus
$\Delta_{(2)}(X_i)$ is obtained and found proportional to $z^2$.
As eq. (\ref{O2}) can be easily generalized to
 \be\label{Om}
 {\cal O}_{(m)}(X_i)=
\Delta_{(m)}(X_i)+ {\cal O}_{(m+1)}(X_i) , \ee \noindent we have
now \be\label{coproduct3} \Delta(X_i)= \Delta_{(0)}(X_i)+ z\;
\delta (X_i) +  \Delta_{(2)}(X_i) + \Delta_{(3)}(X_i) + {\cal
O}_{(4)}(X_i), \ee \noindent where $\Delta_{(2)}(X_i)$ is known
and $\Delta_{(3)}(X_i)$ must be found solving eq.
(\ref{recursionk}) for $k=3$. After a new elimination of unwanted
contributions, a $z^3$-proportional $\Delta_{(3)}(X_i)$ is thus
obtained.

The procedure can be thus iterated obtaining --after that all
$\Delta_{(l)}(X_i)$ for $l<m$ have been obtained in the same way--
$\Delta_{(m)}(X_i)$ that is found to be a polynomial of degree
$m+1$ in $X_j\otimes 1$ and $1\otimes X_j$ and proportional to
$z^m$.

Once all the $\Delta_{(m)}$ are known, the order-by-order
commutation relations are obtained from the homomorphism relation
(\ref{homomorphismproperty}) and, finally, the full coproducts and
commutation relations are obtained as  formal series.

\sect{Standard   quantization of $su(2)$}
\label{quantizationofsu(2)}

To enlighten the details of the construction we discuss explicitly
the standard deformation of $su(2)$. 
The standard ($su(2),\delta $) bialgebra, in the Cartan basis
$\{H,X,Y\}$, is given by the cocommutator map
($\delta : g \to g\otimes g$)
\[
\delta (H)=0 ,\qquad \delta (X)=  H \wedge X ,
\qquad \delta (Y)=   H \wedge Y ,
\]
and the commutation rules \be \label{cartansu2} [H,X]= X, \qquad
[H,Y]= -Y, \qquad [X,Y]= 2\,H. \ee

As stated in Sect.2, we begin to search of the coproducts,
starting from the Lie coalgebra $\Delta_{(0)}$ and finding the
$\Delta_{(k)}$ imposed by $\delta$ to be consistent with
eq. (\ref{recursionk}).

The case of $H$ is simple: we start with $\Delta_{(0)}(H) = H
\otimes 1 + 1\otimes H$ and $\delta(H)=0$, that implies that the
anti-cocommutative part of $\Delta_{(1)}(H)$ is zero. The eq.
(\ref{recursionk}) for $k=1$ (see
Sect.~\ref{Friedrichstheoremrevisited} for details) gives the solution, 
\begin{eqnarray}
\Delta_{(1)}(H)=\alpha_1\,H\otimes H + 
\alpha_2\,(H\otimes X + X\otimes H) + 
\alpha_3\,(H\otimes Y + Y\otimes H) + \nonumber\\
\alpha_4\,X\otimes X +
\alpha_5\,(X\otimes Y + Y \otimes X) +
\alpha_6\,Y\otimes Y
\end{eqnarray}
Since these cocommutative contributions are not 
related to $\delta(H)$ (that in this case vanishes) and the Hopf algebra axioms are fulfilled whatever the $\alpha_i$ coefficients are, the
analytical approach implies that $\alpha_i=0$.
Thus, in agreement with formula
(\ref{coproductseriesdelta}), we write
 \[
 \Delta(H)=
\Delta_{(0)}(H) + {\cal O}_{(2)}(H), \] or, from eq.(\ref{Om}),
\[\Delta(H)=
\Delta_{(0)}(H) + \Delta_{(2)}(H)+ {\cal O}_{(3)}(H),
 \]
and as eq.
(\ref{recursionk}) for $k=2$ is also consistent with
$\Delta_{(2)}(H)=0$ we have
 \[\Delta(H)=
\Delta_{(0)}(H) + \Delta_{(3)}(H)+ {\cal O}_{(4)}(H), \] where the
procedure can be repeated. Thus,  for all the  orders, the
analytical prescription  imposes  $\Delta_{(k)})(H)=0,\;\;
\forall k>0$. Hence \be\label{DE}\Delta(H) =\Delta_{(0)}(H) = H
\otimes 1 + 1\otimes H , \ee i.e. to a null $\delta$, the
analytical procedure associates an object with primitive
coproduct. Note that formula
(\ref{DE}) is not, like in \cite{CP}, a possible choice but the
only coproduct consistent with the analytical prescription.

Equivalently, for $\Delta(X)$ we have
\[
\Delta(X)= \Delta_{(0)}(X)+\Delta_{(1)}(X) + {\cal O}_{(2)}(X) .
\]
The coassociativity condition (\ref{recursionk}) for $k=1$ gives 
\begin{eqnarray}
&& \Delta_{(1)}(X) = z \,\delta(X) + 
\beta_1\,H\otimes H + 
\beta_2\,(H\otimes X + X\otimes H) + \nonumber\\ 
&& \quad \beta_3\,(H\otimes Y + Y\otimes H) +
\beta_4\,X\otimes X +
\beta_5\,(X\otimes Y + Y \otimes X) +
\beta_6\,Y\otimes Y
\end{eqnarray}
where $\beta_i$ are arbitrary constants which are by no means related to $\delta$. As discussed before (see also 
Section 6), we put $\beta_i=0$ and we have
\[
\Delta(X)= \Delta_{(0)}(X)+z \,\delta(X) + \Delta_{(2)}(X_i)+{\cal
O}_{(3)}(X) .
\]
The coassociativity condition (\ref{recursionk}) for $k=2$ solved
in the unknown $\Delta_{(2)}(X)$ gives (again disregarding
arbitrary cocommutative contributions independent from $\delta$):
 \be\label{D2X}
 \Delta_{(2)}(X)= \frac{z^2}{2} (H^2 \otimes X + X \otimes H^2) .
 \ee
By repeating this machinery, we write
\[
\Delta(X)= \Delta_{(0)}(X)+z \,\delta(X) + \Delta_{(2)}(X_i)+
\Delta_{(3)}(X_i)+{\cal O}_{(4)}(X) .
\]
where now $\Delta_{(2)}(X_i)$ is given by eq. (\ref{D2X}) and
$\Delta_{(3)}$ is the new unknown. The coassociativity condition
(\ref{recursionk}) for $k=3$ gives
 \be\label{D3X}
 \Delta_{(3)}(X)= \frac{z^3}{6} (H^3 \otimes X - X \otimes H^3)
 \ee
and the general formula is obtained by iteration, by neglecting order by order the cocommutative contributions unrelated to $\delta$:
 \[
 \Delta_{(k)}(X)= \frac{z^k}{k!} (H^k \otimes X +(-1)^k  X \otimes H^k) \qquad
 \forall k .
 \]
\noindent Now, the $\Delta_{(k)}$ are easily summed to
\[
 \Delta(X)= e^{z H} \otimes X + X \otimes e^{-z H} .
 \]

The approach is exactly the same for $Y$ and gives a similar
result. Thus, we obtain the analytical quantum coproduct
associated to $(su(2),\delta)$
 \be
\label{su2cuantizacion1}\begin{array}{l}
 \Delta (H) = H \otimes 1 + 1 \otimes H\\[0.3cm]
 \Delta (X) = e^{z\,H} \otimes X+ X \otimes e^{-z\,H} \\[0.3cm]
 \Delta (Y) = e^{z\,H} \otimes Y+ Y \otimes e^{-z\,H} ,
\end{array}\ee
that, by inspection, are invariant under the combination of flip
and $z \to -z$.

 Now we have simply to  start from the
commutators (\ref{cartansu2}) and to impose order by order the
homomorphism condition for the deformed commutation rules. The
quantum commutation rules \be \label{su2cuantizacion2} [H,X]= X ,\qquad [H,Y]=
-Y,\ee are quite easy to find. 
The remaining one reads
\be\label{su2cuantizacion3} [X,Y]=\frac{1}{z}\,\sinh(2\,z\,H) ,
 \ee 
 which is a combined result of eqs.~(\ref{su2cuantizacion2}) and
(\ref{homomorphismproperty}). Note that 
 for
$k=0$ the equation (\ref{su2cuantizacion3}) has to give $[X,Y]= 2\,H$. This forbids other $z$-dependent
commutation rules like
\be\label{qquantization}
 [X,Y]=\frac{\sinh(2\,z\,H)}{\sinh z} .
 \ee

Eqs. (\ref{su2cuantizacion1}), (\ref{su2cuantizacion2}) and
(\ref{su2cuantizacion3}) define uniquely the analytical
deformation of the Cartan basis of $su(2)$ such that the
$q$-generators $H,X,Y$ could be called the $q$-Cartan basis of
$su_q(2)$.  By inspection, all the symmetries (for example, $\{H,X,Y\}
\leftrightarrow \{H,-X,-Y\}$) and the embedding conditions  (for
instance, $su(2) \supset \mbox{borel}(H,X) \supset u(1)$) of the
bialgebra $(g,\delta)$ are automatically preserved in the
quantization $(g_q,\Delta)$.

The results of this Section
show that, among all possible coproducts, analyticity chooses the
only one invariant under the combination of flip and change of
sign in $z$ \cite{ballesteros04}.  
Moreover, a $q$-Cartan basis is determined by (\ref{su2cuantizacion1}), (\ref{su2cuantizacion2}) and
(\ref{su2cuantizacion3}), in contradistinction to the usual commutation rule (\ref{qquantization}).

Notice that the coalgebra (\ref{su2cuantizacion1}) is
consistent with other Lie limits as, for instance, $E(2)$ or, after relabeling the generators, with the
twisted jordanian deformation $su_h(2)$ \cite{ohn}. It is only
when also the commutators are included in the game that the
one-to-one correspondence between the bialgebra and the quantum
algebra is obtained.

As a result, we have obtained a Hopf algebra in which 
the coproduct map is such that $z\,\Delta$ is a function of $z\,X_j\otimes 1$ and  $1\otimes z\,X_j$,
and the commutation rules fulfill that $z[X_l,X_m]$ is a function of $z\, X_j$. This is a general property  
of the analytical quantization.

\sect{Standard   quantization of $u(3)$}
\label{quantizationofu(3)}

In the previous section we have described the simplest case of
$su_q(2)$. The procedure above described can be applied to any
bialgebra. As a true example we give now the standard deformation
of $u(3)$. For simple Lie algebras the usual description is made
in terms of the Cartan subalgebra, simple roots and the $q$-Serre
relations without any reference to non-simple roots that remain
undefined \cite{fuchs}. This is a problem for applications where
simple and non-simple roots play the same role. We start instead
from the Weyl-Drinfeld basis of the bialgebra where all roots are
well defined \cite{ballesteros06,ballesteros07} and we obtain a
complete description of the whole structure for $u_q(3) \equiv
su_q(3)\oplus u(1)$, real form of $A_2^q\oplus A_1$. In this
basis, the explicit commutation rules are ($i,j,k=1,2,3$):
\[\begin{array}{l}\label{commutatorsA}
[H_i,H_j]=0, \\[0.3cm]
[H_i,F_{jk}]=(\delta_{ij} - \delta_{ik})F_{jk}, \\[0.3cm]
[F_{ij},F_{kl}]=(\delta_{jk} F_{il}- \delta_{il} F_{kj}) +
\delta_{jk} \delta_{il}(H_i - H_j) .
\end{array}\]

The canonical Lie bialgebra structure  is determined by the
cocommutator: \be\begin{array}{l}\label{cocommutadorA}
\delta (H_i)=0,\\[0.3cm]
\delta (F_{ij})= \frac{z}{2}\; (H_i-H_j) \wedge F_{ij} +z\,
\sum_{k=i+1}^{j-1}{F_{ik} \wedge F_{kj}}
\qquad\qquad (i<j),\\[0.3cm]
\delta (F_{ij})=  \frac{z}{2}\; (H_j-H_i) \wedge F_{ij} - z\,
\sum_{k=j+1}^{i-1}{F_{ik} \wedge F_{kj}}   \qquad\qquad (i>j).
\end{array}\ee

We begin with the coalgebra of the Borel subalgebra $b_+ \equiv \{
H_1, H_2, H_3, F_{12}, F_{13}, F_{23}\}$. Repeating the procedure
of the preceding paragraph (or, simply, remembering the embeddings
$su_q(3) \supset su_q(2)$) we get
\be\begin{array}{l}\label{cocommutadorB}
 \Delta (H_i) = H_i \otimes 1 + 1 \otimes H_i ,\\[0.3cm]
 \Delta (F_{12}) = e^{\z\,(H_1-H_2)/2} \otimes F_{12} + F_{12} \otimes
e^{-\z\,(H_1-H_2)/2} ,\\[0.3cm]
\Delta (F_{23}) = e^{\z\,(H_2-H_3)/2} \otimes F_{23} + F_{23}
\otimes e^{-\z\,(H_2-H_3)/2} \,.
\end{array}\ee
\noindent The explicit quantization of $\Delta(F_{13})$ from
eq.(\ref{cocommutadorA}) requires more work. We find
\be\begin{array}{l}\label{DF13}
 \Delta(F_{13})= \; \,
 e^{\z\,(H_1-H_3)/2} \otimes F_{13} + F_{13} \otimes
e^{-\z\,(H_1-H_3)/2} \; \; \; \; + \hskip0.2cm \; 2 \,\sinh \frac z2 \;
\times\,\\[0.3cm]
\displaystyle \left( e^{\z\,(H_2-H_3)/2}\,F_{12}
 \otimes e^{-\z\,(H_1-H_2)/2}\,F_{23} -
\, e^{\z\,(H_1-H_2)/2}\,F_{23} \otimes
e^{-\z\,(H_2-H_3)/2}\,F_{12}\right) ,
\end{array}\ee
\noindent that (like \ref{cocommutadorA}) is inconsistent with
usual definition \cite{alisauskas-smirnov} \be\label{F13}
F_{13}':=e^{z/2}\,F_{12}\,F_{23} - e^{-\z/2}\,F_{23}\,F_{12} . \ee
\noindent Cartan matrix, $q$-Serre relations and $q$-commutators
do not seem perhaps the simplest approach to quantum algebras. The
origin of the definition (\ref{F13}) is indeed related to the
$q$-Serre relations \be\begin{array}{l}\label{qserre}
 F_{12}^2\,F_{23} - (e^{\z}+e^{-\z})\,F_{12}\,F_{23}\,F_{12}
+ F_{23}\,F_{12}^2=0 ,\\[0.3cm]
 F_{12}\,F_{23}^2 - (e^{\z}+e^{-\z})\,F_{23}\,F_{12}\,F_{23}
+ F_{23}^2\,F_{12}=0 ,
\end{array}\ee
that, using (\ref{F13}), can be written
\[\begin{array}{l}
 e^{\z/2}\,F_{23}\,F_{13}' - e^{-\z/2}\,F_{13}'\,F_{23} =0 ,\\[0.3cm]
 e^{\z/2}\,F_{13}'\,F_{12} - e^{-\z/2}\,F_{12}\,F_{13}' =0 ,
\end{array}
\]
showing that $F_{13}'$ $q$-commutes with both $F_{12}$ and
$F_{23}$. Anyway, imposing the homomorphism
eq.(\ref{homomorphismproperty1}) of the coproducts, we find that
the commutators
 \[
  [F_{12},F_{23}] = F_{13}, \qquad [F_{32},F_{21}] =
F_{31}
\]
\noindent of the bialgebra remain unchanged in the quantization.
In agreement with the quantum theories, where the commutator is
connected to the measure, the commutator remains the appropriate
map also  in the deformation of $g$.

The quantized coproduct of $b_- \equiv \{ H_1, H_2, H_3, F_{21},
F_{31}, F_{32}\}$ is similar:
 \[\label{bm1}\begin{array}{l}
 \Delta (H_i) = H_i \otimes 1 + 1 \otimes H_i ,\\[0.25cm]
 \Delta (F_{21}) = e^{\z\,(H_1-H_2)/2} \otimes F_{21} + F_{21} \otimes
e^{-\z\,(H_1-H_2)/2}  ,\\[0.25cm]
 \Delta (F_{32}) = e^{\z\,(H_2-H_3)/2} \otimes F_{32} + F_{32} \otimes
e^{-\z\,(H_2-H_3)/2} ,
\end{array}\]
\[\label{bm2}\begin{array}{l}
 \Delta(F_{31}) = \; \; \;
 e^{\z\,(H_1-H_3)/2} \otimes F_{31} + F_{31} \otimes
e^{-\z\,(H_1-H_3)/2} \; \; \; + \\[0.3cm]
\; \displaystyle 2\, \sinh \frac z2 \,\left(\,
e^{\z\,(H_2-H_3)/2}\,F_{21}
 \otimes e^{-\z\,(H_1-H_2)/2}\,F_{32} -
\, e^{\z\,(H_1-H_2)/2}\,F_{32} \otimes
e^{-\z\,(H_2-H_3)/2}\,F_{21}\right) .
\end{array}\]
\noindent All coalgebra is thus known. Now we have to find the
deformed commutation rules compatible with the above coalgebra and
the $u(3)$ limit.

From $q$-Serre relations (\ref{qserre}) (as well as from the
$\Delta$ isomorphism) we get the only commutation rules for
$b_\pm^q$ that are deformed:
\[\begin{array}{l}
 \conm{F_{12}}{F_{13}}=\conm{F_{12}}{\conm{F_{12}}{F_{23}}}=4\,
 (\sinh \frac z2)^2\, F_{12}\,F_{23}\,F_{12},\\[0.25cm]
 \conm{F_{13}}{F_{23}}=\conm{\conm{F_{12}}{F_{23}}}{F_{23}}=4\,
 (\sinh \frac z2)^2 \,  F_{23}\,F_{12}\,F_{23} ,\\[0.25cm]
 \conm{F_{31}}{F_{21}}=
\conm{\conm{F_{32}}{F_{21}}}{F_{21}}= 4\,
 (\sinh \frac z2)^2\, F_{21}\,F_{32}\,F_{21},\\[0.25cm]
 \conm{F_{32}}{F_{31}}=
\conm{F_{32}}{\conm{F_{32}}{F_{21}}}= 4\,
 (\sinh \frac z2)^2 \,  F_{32}\,F_{21}\,F_{32} .
\end{array}\]
Now we have to consider the crossed commutation relations. We
start from $\conm{F_{23}}{F_{21}}$. The compatibility with the
quantum coproduct leads to the equation
\[
\Delta(\conm{F_{23}}{F_{21}})= e^{\z\,(H_1-H_3)} \otimes
\conm{F_{23}}{F_{21}} + \conm{F_{23}}{F_{21}} \otimes
e^{-\z\,(H_1-H_3)} .\] Thus, in agreement with \cite{fuchs}, the
unique analytical solution consistent with the coproduct map for
$\conm{F_{23}}{F_{21}}$ and $\conm{F_{12}}{F_{32}}$ are:
 \[
 \conm{F_{23}}{F_{21}}=0 ,\qquad
\conm{F_{12}}{F_{32}}=0.
\]
From the two embedded $su_q(2)$ Hopf subalgebras we get
\[
[F_{12},F_{21}]= \frac{1}{z}\, \sinh(\z\,(H_1-H_2)),\qquad
[F_{23},F_{32}]=\frac{1}{z}\,\sinh(\z\,(H_2-H_3)),
 \]
and, because the full structure written in terms of the
commutation rules, the (not deformed) Jacobi identities can be
used as a short cut to derive
\[\begin{array}{l}
 [F_{13},F_{21}]=\conm{\conm{F_{12}}{F_{23}}}{F_{21}}=
 -\conm{\conm{F_{21}}{F_{12}}}{F_{23}}=\\[0.25cm]
\hskip0.5cm \frac{1}{z} \conm{\sinh(\z\,(H_1-H_2))}{F_{23}}=
-\frac{2}{z}
\,\sinh\frac{z}{2}\,\cosh(\z\,(H_1-H_2+\frac12))\,\,F_{23} ,
\end{array}\]
and analogously
\[\begin{array}{l}
[F_{13},F_{32}]=\frac{2}{z} \,\sinh\frac{z}{2}\,\cosh(\z\,
(H_2-H_3+\frac12))\,F_{12},\\[0.25cm]
[F_{12},F_{31}]= -\frac{2}{z} \,\sinh\frac{z}{2}
\,\cosh(\z\,(H_1-H_2
-\frac12))\,F_{32},\\[0.25cm]
[F_{23},F_{31}]=\frac{2}{z} \,\sinh\frac{z}{2}\,\cosh(\z\,
(H_2-H_3-\frac12))\,F_{21}.
\end{array}\]
The last relation is computed imposing the homomorphism property,
obtaining
 \[\begin{array}{l}
 [F_{13},F_{31}]= \frac{1}{z} \sinh(\z\,(H_1-H_3))
 +\frac{2}{z} (\sinh\frac{z}{2})^2 \,\sinh(\z\,(H_1-H_2))\,
 \anticonm{F_{23}}{F_{32}}+ \\[0.15cm]
\hskip 2.5cm\frac{2}{z}
(\sinh\frac{z}{2})^2\,\sinh(\z\,(H_2-H_3))\,
\anticonm{F_{12}}{F_{21}} .
\end{array}\]

\sect{Friedrichs theorem revisited: $U(g) \rightarrow
(g,\Delta_{(0)})$}\label{Friedrichstheoremrevisited}

The universal enveloping algebra $U(g)$ is defined in terms of an arbitrary set of $n$ basic
elements $\{Y_j\}$ on which a PWB basis for the whole $U(g)$ can
be built. They are not (in principle) primitive but they are
cocommutative.

Here we give a constructive proof of the Friedrichs theorem,  building
explicitly the primitive generators  $\{X_j\}$ in terms of the
$\{Y_j\}$. The machinery consists in repeated changes of bases
that allow to obtain each time  a better approximation to
primitivity where the problem is reformulated at each step in
terms of the preceding basis. An infinite iteration of the
procedure allows to find,  among the $\infty$-many possible bases
in the $U(g)$, the Lie generators.

More explicitly,  we consider that $X_i \equiv \lim_{k\to \infty} X_i^k $ were
$\{X_i^k\}$ is a basic set that approximates the Lie-Hopf
coproducts up to ${\cal O}_{(k)}(X_i)$. The   terms ${\cal
O}_{(m)}(Z_i)$, defined in Sec.~\ref{analyticalquantization}  as
series of degree greater that $m$ in $Z_j\otimes 1$ and $1\otimes Z_j$, are,
in this Section, cocommutative  since we are working in $U(g)$.

To begin with, let us define an homogeneous symmetric polynomial
of order $m$:
\[ P_{(m)}(Z_i) :=\sum_{m_i}
f_i^{m_1,m_2\dots m_n} \;{\cal S}_{m} \left[
(Z_1)^{m_1}(Z_2)^{m_2}\cdots (Z_{n})^{m_n} \right]  ,
\]
where the
sum on the $m_i$ is restricted to $\sum m_i=m$, $f_i^{m_1,m_2\dots
m_n} \in \C$ and
\[
 {\cal S}_m \left[(Z_1)^{m_1}(Z_2)^{m_2}\cdots
(Z_{n})^{m_n}\right]\equiv \sum_{\sigma \in {\bold S}_m}
\sigma\left[(Z_1)^{m_1}(Z_2)^{m_2}\cdots (Z_{n})^{m_n}\right] ,
\]
being ${\bold S}_m$   the group of permutations of order $m$.

Now any original basic set $\{Y_i\}$ is a zero approximation to
$\{X_i\}$ : $X_i^0 := Y_i$. Indeed
\be\label{aa}
 \Delta(X_i^0)
= \Delta_{(0)}(X_i^0) + {\cal O}_{(1)}(X_i^0) = X_i^0\otimes 1 +
1\otimes X_i^0 + {\cal O}_{(1)}(X_i^0).
\ee
The explicit form of
${\cal O}_{(1)}(X_i^0)$ in (\ref{aa}) is
\[\label{O11}
  {\cal O}_{(1)}(X_i^0)
= \sum c_i^{jl}\; \left( X_j^0 X_l^0\otimes 1
+ 1\otimes X_j^0 X_l^0 \right)+\sum d_i^{jl}\; X_j^0\otimes X_l^0
 \;+ \;{\cal O}_{(2)}(X_i^0)
\]
where $c_i^{jl}$ and $d_i^{jl} = d_i^{lj}$ are constants. Again
from (\ref{coproductseries}) , \be\label{O1}
  {\cal O}_{(1)}(X_i^0)= \Delta_{(1)}(X_i^0)+
{\cal O}_{(2)}(X_i^0),\ee
 and we have to impose on $\Delta_{(1)}(X_i^0)$
the coassociativity condition (\ref{recursionk}) for $k=1$ that
gives $c_i^{jl}\equiv 0$ , while no more restrictions are found on
$d_i^{jl}$.
\noindent Thus, if we define
 \be\label{vv} P_{(2)}(X_i^0) :=
\sum d_i^{jl}\; [X_j^0, X_l^0]_{+} \ee
\noindent we have \be\label{B} \Delta_{(1)}(X_i^0) =
 \Delta_{(0)}(P_{(2)}(X_i^0))-P_{(2)}(X_i^0)\otimes 1-1\otimes P_{(2)}(X_i^0) .\ee
\noindent We can thus define the next approximation of the Lie
generators \be X_i^1 := X_i^0 -\; P_{(2)}(X_i^0)\label{0to1} \ee
\noindent and we get for $\{X_i^1\}$ a coproduct with {\it
vanishing first order contributions}: \[
 \Delta(X_i^1) =
\Delta_{(0)}(X_i^1)+ {\cal O}_{(2)}(X_i^0). \]

Still more relevant, eq. (\ref{0to1}) allows to write  ${\cal
O}_{(2)}(X_i^0)$   in terms of
$\{X_i^1\}$ as ${\cal O}_{(2)}(X_i^1)$:
 \[ \Delta(X_i^1) = \Delta_{(0)}(X_i^1) + {\cal O}_{(2)}(X_i^1). \]

As both the relation (\ref{O1}) and (\ref{B}) can be generalized
to \[\label{Omm} {\cal O}_{(m)}(X_i^{m-1})=
\Delta_{(m)}(X_i^{m-1})+ {\cal O}_{(m+1)}(X_i^{m-1}),\]
\[\label{BB} \Delta_{(m)}(X_i^{m-1}) =
 \Delta_{(0)}(P_{(m+1)}(X_i^{m-1}))-P_{(m+1)}(X_i^{m-1})\otimes 1
 -1\otimes P_{(m+1)}(X_i^{m-1}) ,\]
\noindent we are ready for next step. Imposing the coassociativity
property on the most general symmetric polynomial of third order
in \,$X_j^{1}\otimes 1$\, and \,$1\otimes X_j^{1}$\,, we get
\be\label{CC} {\cal O}_{(2)}(X_i^1)= \Delta_{(2)}(X_i^1)+ {\cal
O}_{(3)}(X_i^1),\ee
 \be\label{DD}
\Delta_{(2)}(X_i^1)=\Delta_{(0)}(P_{(3)}(X_i^1))-P_{(3)}(X_i^1)\otimes
1-1\otimes P_{(3)}(X_i^1) , \ee

With a new change of basis $ X_i^2 := X_i^1 - P_{(3)}(X_i^1) $ we
obtain the coproduct of the second approximation $\{X_i^2\}$ to
the generators in terms of the same $\{X_i^2\}$
\[
\Delta(X_i^2) = \Delta_{(0)}(X_i^2) + {\cal O}_{(3)}(X_i^2) ,
\]
now {\it free form both first and second order contributions}.

The procedure can now be easily iterated and the
$\Delta_{m}(X_i^{m-1})$ contribution eliminated through a new
change of basis that affects the higher orders only. The residual
term becomes $O_{(m+1)}(X_i^m)$ and we get the $m$-order
approximation to the Lie generators
\[
 \Delta(X_i^m) = \Delta_{(0)}(X_i^m) + {\cal O}_{(m+1)}(X_i^m).
\]

The true generators of the Lie algebra $g$ are (formally)
recovered in the limit \[ X_i := \lim_{m\to \infty} X_i^m \] and,
in agreement with the Friedrichs theorem, their coproduct is the
primitive one \[ \lim_{m\to \infty} \Delta(X_i^m) =\lim_{m\to
\infty} \Delta_{(0)}(X_i^m) =\Delta_{(0)}(X_i) = \Delta(X_i) =
X_i\otimes 1 + 1\otimes X_i. \]
 Of course, this coproduct is an
algebra homomorphism with respect to the (linear) Lie algebra
commutation rules, \[ \Delta[X_i, X_j] = [X_i,X_j]\otimes 1 +
1\otimes [X_i,X_j] \] and the $n$ Lie algebra generators are
univocally identified in a constructive manner within $U(g)$
pushing away order by order the corrections to a primitive
coproduct.

Let us stress again that the central point of this analytical
approach to Friedrichs theorem (as well as to its following
extension to $U_q(g)$) is that at each order all relations can be
rewritten in terms of the corresponding approximations of the
generators.


\sect{Extension of Friedrichs theorem: $U_q(g) \rightarrow
(g_q,\Delta)$}\label{ExtensionFriedrichstheorem}

Of course, Friedrichs theorem is a well known result. However, it has been described here
because the procedure that allows to individuate the quantum
algebra generators $g_q$, among the $\infty$-many possible bases
of the $U_q(g)$ is exactly the same that allows to individuate the
generators $g$ among the $\infty$-many possible bases of the
$U(g)$: therefore, the analytical approach can be considered as an extension
to quantum algebras of the Friedrichs theorem.

The preceding construction indeed works also in the $\delta\neq 0$
case, thus providing us the prescription for the construction of
the almost primitive generators --obtained in
Sec.~\ref{analyticalquantization} from $(g, \Delta_{(0)})$--
starting from an arbitrary set of basic elements of any $U_q(g)$.

As before, we start with an arbitrary set of basic elements
$\{Y_j\}$ that define the $U_q(g)$ (and no more a $U(g)$) with, as
classical limit, a Lie bialgebra with $\delta \neq 0$.

Eqs. (\ref{aa}) and (\ref{O1}) are still valid but now to the
$\Delta_{(1)}$ of (\ref{B}) we have to add the contribution of
$\delta$. We have thus
\[ \Delta(X_i^0)= \Delta_{(0)}(X_i^0)+ z \,\delta(X_i^0) + \sum
d_i^{jl}\; X_j^0\otimes X_l^0 \; + {\cal O}_{(2)}(X_i^0). \]

The same $P_{(2)}(X_i^0)$ of eq. (\ref{vv}) allows us to define
again $X_i^1 := X_i^0 - P_{(2)}(X_i^0)$. As this change of
variables does effect the $\delta$ contribution only to higher
orders, the differences can be included in ${\cal O}_{(2)}(X_i^0)$
(or, equivalently, ${\cal O}_{(2)}(X_i^1)$). We can thus write
\[ \Delta(X_i^1) = \Delta_{(0)}(X_i^1) + z\, \delta(X_i^1) +
{\cal O}_{(2)}(X_i^1). \] Next step is, like in eq. (\ref{CC}),
 to introduce
$\Delta_{(2)}(X_i^{1})$
 \[
{\cal O}_{(2)}(X_i^{1})= \Delta_{(2)}(X_i^{1})+ {\cal
O}_{(3)}(X_i^{1}). \] 
Like in the $\Delta_{(1)}$ case, $\Delta_{(2)}$ has two contributions:
one of them proportional
to $z^2$ (and imposed by the consistency between $\delta$ and
coassociativity) and the other one described in (\ref{DD}). Again, the latter can be removed by another change of basis that does not modify the
form of the $z$-dependent contributions as the induced modifications can
be included in ${\cal O}_{(3)}(X_i^2)$.

This procedure can be iterated. Once the problem is solved for
$\Delta_{(m-1)}$, a contribution to  $\Delta_{(m)}$, proportional
to $z^m$ (cocommutative for $m$ even and anti-cocommutative for
$m$ odd), is found while the unessential $z$-independent terms are
removed exactly as in the case $\delta = 0$, with a change of
basis that does not affect the form of known $z$-depending terms
because the introduced changes are always to orders higher of
$z^m$ and thus pushed out in ${\cal O}_{(m+1)}$. For $m \to
\infty$ the same almost primitive coproducts derived from
$(g,\Delta_{(0)})$ in Sect.~\ref{analyticalquantization} are found
to be a basic set for $U_q(g)$. The homomorphism condition imposes, of course,
the same deformed commutation rules and we have thus closed the
other side of the diagram, finding that $U_q(g)$ has as one of its
basic sets the same almost primitive set $(g_q,\Delta)$ that we have obtained
before by analytical continuation of the Lie generators.

Note that the deformation does not affect the fundamental trick of
the game: the iterative procedure where each order is stated in
terms of the preceding ones.

As we have explicitly constructed the transformation between any
arbitrary basic set $\{Y_j\}$ and the quantum algebra generators
$\{X_j\}$, we have demonstrated the one-to-one correspondence
between $U_q(g)$ and $(g_q,\Delta)$.
Since in Sect.~\ref{analyticalquantization} it has been
shown that $(g_q,\Delta)$ is in one-to-one correspondence with
$(g,\delta)$, the classification of quantum groups has been
reduced to the classification of Lie biagebras.


\section{Concluding remarks}

The main result of this paper is the construction of the unique
almost primitive  basic set  characterizing the quantum universal
enveloping algebra $U_q(g)$, in perfect analogy with the unique
primitive basic set determining $U(g)$. As the last one is the Lie-Hopf
algebra $(g,\Delta_{(0)})$, we call the first one {\it quantum
algebra} $(g_q,\Delta)$.  Hence, a  deformed structure of the same
dimension that the underlying Lie algebra is introduced instead of
the $\infty$-dimensional quantum universal enveloping algebra.
This quantum algebra could be the essential object  to be connected with
physical operators, like in the Lie case where the generators do
have a precise meaning in terms of symmetry transformations.

We have also shown
that the connections between bialgebras, quantum algebras and
quantum universal enveloping algebras are always one-to-one, such
that the classification problem (as well as equivalence relations and
embedding properties) can be stated at the Lie bialgebra level.
As a third point, note also that the analytical quantization method here presented is constructive and
it could be implemented  by making use of computer algebra. 

Besides the proposed almost primitive basis, we would like to quote two other
relevant bases that play a role both in mathematics and in
physics: the Lie basis and the canonical/crystal basis. In the Lie
basis (for instance, in twisting) the algebra remains unmodified
and all deformations affect the coalgebra offering a possible way
to introduce an interaction but saving the global symmetry
\cite{enrico03}. In the canonical or crystal basis (with
applications in statistical mechanics \cite{Kang} and in genetics
\cite{MS})), instead, the algebraic sector of the Hopf algebra
structure is obtained in the limit $|z| \to \infty$ \cite{Lu,Ka}
and coalgebra is a byproduct.

We stress that, within this approach, the fundamental object for the construction of the quantum algebra
is the coproduct, while the deformed
commutation rules are derived {\it a posteriori} by making use of
the homomorphism property. It is also worthy noting that in the usual quantization of simple
Lie algebras --as the whole structure is defined in terms of the
Cartan subalgebra and simple roots-- the $q$-generators associated
to non-simple roots (that, for a physics are as relevant as the
others) do not play any role. As a consequence, these
$q$-generators can be defined in many different ways, in
contradistinction to the Lie case. However the analyticity, here
introduced, forbids the $q$-commutator as the appropriate bracket
in $U_q(g)$, and non-simple root generators are uniquely defined
--exactly as in the Lie case-- by commutation relations.

The fact that the quantum algebra is built up in terms of
commutators makes possible a straightforward semiclassical limit in terms of
Poisson-Lie structures. In particular,
applications of the Poisson $su_q(3)$ algebra will be presented elsewhere.


\section*{Acknowledgments}

This work was partially supported  by the Ministerio de
Educaci\'on y Ciencia  of Spain (Projects FIS2005-03989 and FIS2004-07913),  by the
Junta de Castilla y Le\'on   (Project VA013C05), and by
INFN-CICyT (Italy-Spain).




\end{document}